\documentclass[11pt]{amsart}
\usepackage{amsmath,amssymb,amscd}
\usepackage{amsmath,amssymb}
\usepackage{mathrsfs}
\usepackage{color}
\usepackage{amsxtra, amsmath}
\usepackage{amssymb, amscd}
\usepackage{graphicx}

\newcommand\Erf{\textrm{Erf}}

\theoremstyle{plain}

\theoremstyle{definition}

\theoremstyle{remark}

\title[The bispectral problem, the Darboux process, monodromy and the Hermite operator]{The bispectral problem, the Darboux process, monodromy and the Hermite operator}

\author{M. M. Castro \and  F. A.  Gr\"unbaum}

\address{ Departamento de Matem\'atica Aplicada II and IMUS, Escuela Politécnica Superior, Universidad de Sevilla, 41011, Sevilla, Spain}
\email{mirta@us.es}
\thanks{The research of M. M. Castro was partially supported by  PID2024-155593NB-C21  (FEDER(EU) / Ministerio de Ciencia, Innovaci\'on y Universidades-Agencia Estatal de Investigaci\'on) and FQM-262 (Junta de Andalucía).
	 }

\address{Department of Mathematics, University of California, Berkeley
CA 94705}
\email{grunbaum@math.berkeley.edu}



\begin{document}

\begin{abstract} The complete solution of the bispectral problem for the Schr\"odinger operator $L=-\tfrac{d^2}{dx^2}+V(x)$ in \cite{DG} is obtained by the application of the Darboux process to the cases of  $V=0$ and $V(x)=-\tfrac{1}{4x^2}$. Both of these cases are trivially bispectral and after repeated applications of the Darboux process one gets either a pair of rank one bundles of bispectral situations (when starting from $V=0$) or a rank two bispectral bundle (when starting from $V(x)=-\tfrac{1}{4x^2}$).  In the first case all operators have ``trivial monodromy'' as defined in \cite{DG}. In the second case the monodromy group of all operators is given by the integers.  In this paper we start from $V(x)=x^2$, use the Darboux process and explore the connection between the rank of certain non-polynomial bispectral families and trivial monodromy by means of examples. The main conclusion is that the results in \cite{DG} do not apply verbatim in this case.


\end{abstract}

\maketitle

\section{A brief  introduction}

	The bispectral property discussed in \cite{DG} can be seen as a relative of a problem discussed, among others, by S. Bochner, \cite{Boch}. For a very nice review of this material see the paper by Luc Haine \cite{Haine1}, where one finds a discussion of an extension studied by H. L. Krall and A. M. Krall (father and son). See also \cite{HK} for a nice collection of papers. The paper by L. Haine mentioned above contains a lot of new material in terms of applying the Darboux method to {\bf double infinite} matrices, i.e., ignoring the boundary conditions that are imposed in the case of orthogonal polynomials connected with {\bf semi infinite} matrices. An interesting early example of this is given in \cite{GruHaine}. 

\bigskip

The results in \cite{DG} feature two different situations: a pair of ``rank one'' bispectral families and a ``rank two'' bispectral family (see Section \ref{Hermiteclasical} for a formal definition of rank in this framework). Both of these are obtained by repeated applications of what is known as the Darboux process to a second order differential operator of the form
$$ -\dfrac{d^2}{dx^2} + V(x).$$ Incidentally this process is
	actually due to Moutard, as one can read in \cite{Dar}. The result is given in Livre IV, Chap IX, No. 408, and the attribution to Moutard is given in Chap II, No. 343, page 53 of \cite{Dar}. The origin of this powerful method may go further back as mentioned in \cite{Gesz}.

\bigskip

The Darboux process can be applied to difference operators as well, as first done in \cite{Matveev1}, and its usefulness in the context of the bispectral problem was exploited in  \cite{GH0,GH1,GH2,GHH3,GHH4,GruHaine}. Of particular interest in these papers is the case of the Hermite polynomials. This is the one case, among the classical orthogonal polynomials, where the support of the orthogonality measure runs over the entire real line making it quite different from the case of the Krall-Jacobi and the Krall-Laguerre polynomials discussed in the first three papers just cited. The application of the Darboux process is greatly simplified when the orthogonality measure has at least one finite end point, and this feature sets the Hermite case apart from the other ones. This is treated in detail in \cite{GHH4}. See also \cite{Deift} for a way to handle a related issue.

\bigskip

 We notice that the bispectral property in \cite{DG} was initially considered as a step in an ambitious program dealing with the issue of ``time-and-band-limiting'' as motivated by medical imaging in the presence of limited and noisy data. The introduction to \cite{DG} briefly mentions this; for further work one can see \cite{G3,G8,G6} as well as the more recent set of papers
  \cite{CGYZ1,CGYZ2,CGYZ3,CGYZ4} together with \cite{ CG5,CG6,CG24,CGPZ, CGZ, GPZ1, GPZ2, GPZ3}.

\bigskip

 We close this introduction by mentioning that the study of the bispectral problem in the ``continuous-continuous'' case
posed and solved in  \cite{DG}, is full of unexpected surprises that could not 
be surmised from its simple formulation. Here is a short list:

\begin{itemize}
\item[a.] The ad-conditions are necessary and sufficient for bispectrality.
\item[b.] As a consequence, the potential needs to be rational,
one can exclude irregular singularities and the
potential needs to decay at infinity. All of this is proved in \cite{DG}.
\item[c.] The connection with integrable systems (KdV)
\item[d.] The relevance of the Darboux process.
\item[e.] (Most important for the present paper) The relation with trivial or 
integer monodromy, and the strict connection with rank one and rank two
bundles respectively.
\end{itemize}

In particular the
observation about trivial monodromy and the explicit form of the 
potential gave origin to a paper by V. Goncharenko and  A.P. Veselov \cite{GonchVeselov},
where they deal with the matrix valued case of the Schroedinger 
equation and where they say that this
is a result that one could imagine was settled in the previous century
but was only done in \cite{DG}.
Reference \cite{Oblom} by A. Oblomkov extends this result
from \cite{DG} and gives the class
of potentials with quadratic growth at infinity that have trivial monodromy.    As indicated in the present paper our examples fall within this class.

The work in \cite{DG} deals with the continuous-continuous case and sets up a 
clear correspondence between trivial monodromy and rank one bundles on the
one hand and integer monodromy and rank two bundles on the other hand.

The case studied in the present paper
is of mixed kind: a potential that explodes at infinity calls for
a discrete spectrum and in the best of cases for difference equations
satisfied by the eigenfunctions.
To the best of our knowledge this situation regarding monodromy has not
been explored in the 40 years since the appearence of \cite{DG}.
What we find here is that the lessons learned earlier do 
not apply in the case of a very important potential that does
not decay at infinity. 

Here we display examples with trivial monodromy going with rank
two bundles of common eigenfunctions giving a bispectral situation.
To the best of our knowledge, this has never been seen before and stands in
sharp contrast to the results in \cite{DG}.

\section{The contents of the paper}

  We will consider a few examples of differential operators obtained from the Hermite differential  operator 
$\mathcal{H}_x$, namely $$
\mathcal{H}_x=\partial_x^2-2x\partial_x,
$$by applications of the Darboux process, and consider the issue of bispectrality in these situations with some emphasis on the role played by trivial monodromy. We are interested in the full two dimensional space of eigenfunctions of the differential operator and not just in polynomial eigenfunctions.

\bigskip

Our results fit with the spirit of \cite{GH1,GH2,Reach123,R} where the issue of bispectrality in the continuous-discrete case was first considered. Many of the tools relevant to this mixed case are related to those  given in \cite{DG};
besides the Darboux process we note the use of
the so called ad-conditions and the role of trivial or very simple monodromy. Some of these tools have proved to be of use in the more
recent study of exceptional orthogonal polynomials, see for instance \cite{GGU2,GU1,GUGrMi,GUKKM,Qu,STZ}. See, in particular, both Section 2.3 and the closing remarks in \cite{GUKKM}. These issues feature in some of the discussions below.

One more comment is appropriate here: the big surprise in  \cite{DG} is the connection between the bispectral problem and integrable systems such as KdV. Many of the tools in this paper were ``floating around'' in connection with developments around the KdV equation around the 1970s, and this allowed the authors to use these tools. There are two tools that appear for the first time in this paper: the relevance of the ad-conditions and of monodromy considerations in connection with the bispectral property, in particular the characterization given in Proposition 3.3 and its consequence Theorem 3.4, page 196. One might imagine that this Proposition is a result from the  19th century, however see \cite{FeldVa, GonchVeselov,Oblom}. Interesting extensions to the matrix valued case have been given by A. Veselov and his collaborators. This result has been used, in the scalar valued case, for instance in \cite{GomGrdMil7,Oblom}.

\bigskip

In Section $3$ we look at the familiar case of the Hermite equation looking for two dimensional families of bispectral situations. For comparison we mention what happens in the simpler case of the Legendre operator.

\bigskip

Section $4$ is devoted to the case of a simple Darboux transformation from the case of Section $3$, and we find a multiparameter collection of these rank two bispectral families.

\bigskip

A further example is given in Section $5$ where the freedom in picking free parameters is larger than in the previous example. All of these examples have trivial monodromy, showing that the situation in \cite{DG} where simple monodromy went only with rank one bispectral families does not hold.

\section{The Hermite operator}\label{Hermiteclasical}

Consider the Hermite operator
\begin{equation}\label{opHermite}
\mathcal{H}_x=\partial_x^2-2x\partial_x
\end{equation}and its eigenfunctions  $\phi_n(x)$ given by
\begin{equation}\label{eigenf}
(\mathcal{H}_x \phi_n)(x)=-2n\phi_n(x),\quad n=0,1,2,\ldots .
\end{equation}


\bigskip


For each $n$ we have a two dimensional space $\mathcal{V}_n$ of solutions to (2). The subspace $\mathcal{V}_0$ is spanned by the functions $1$ and $\Erf( i x)$ where $\Erf$ denotes the well known entire Error function,

\bigskip

$$\Erf(z)=\dfrac{2}{\sqrt{\pi}}\int_0^{z}e^{-t^2}dt .$$

\bigskip


This situation holds also for $n=1,2,3,...$, namely the space $\mathcal{V}_n$ is spanned
by the Hermite polynomial
$$
H_n(x)=(-1)^ne^{x^2}\partial_x^ne^{-x^2}
$$

\bigskip

\noindent and some nonpolynomial function that can be chosen in many different ways.
For our purpose define the
entire functions
$$
\xi_n(x)=\sum_{m=0}^{n}\dfrac{(-x)^m}{2^{n-m}m!\Gamma(1+\frac{n-m}{2})}+
\sum_{j=0}^{\infty}\dfrac{(-x)^{n+2j+1}2^{2j+1}}{(n+2j+1)!\Gamma(\frac{1}{2}-j)},\quad n\geq 0.
$$
\bigskip

It is easy to see that the functions $\xi_n(ix)$ and $\xi_n(-ix)$ are linearly independent solutions of (\ref{eigenf}).

\bigskip

The expression for the Error function above can be given by the standard method used to produce a second solution of a second order differential equation once one has one solution.

\bigskip


There are, of course, other possible choices of a basis for the solutions of  (\ref{eigenf}), such as the pair 

\begin{equation}\label{basehyperg}
  _1F_1\left(-\frac{n}{2};\frac{1}{2};x^2\right)\quad \textrm{and} \quad x\ _1F_1\left(\frac{1}{2}-\frac{n}{2};\frac{3}{2};x^2\right).
\end{equation}
\bigskip

It might be appropriate to observe that the second solution above can be obtained from the first one by the classical procedure mentioned earlier, namely$$
x\ _1F_1\left(\frac{1}{2}-\frac{n}{2};\frac{3}{2};x^2\right)=\ _1F_1\left(-\frac{n}{2};\frac{1}{2};x^2\right)\int_0^x\dfrac{e^{y^2}}{_1F^2_1\left(-\frac{n}{2};\frac{1}{2};y^2\right)}dy.
$$

The function in the denominator in the integral does not vanish at $y=0$ and the integral makes sense for $x$ small enough.

\bigskip


		The functions $\xi_{n}(x)$ can be written in terms of the confluent hypergeometric function as follows:
		

\bigskip

\begin{equation*}
	\xi_{n}(\pm i x)=\dfrac{1}{\sqrt{\pi}n!}\left[\Gamma \left(\frac{n+1}{2}\right) \, _1F_1\left(-\frac{n}{2};\frac{1}{2};x^2\right)\mp 2 i x \Gamma \left(\frac{n}{2}+1\right) \, _1F_1\left(\frac{1}{2}-\frac{n}{2};\frac{3}{2};x^2\right)\right].
\end{equation*}

One has, in particular, that the first solution of the pair above in (\ref{basehyperg}) with $n$ replaced by $2n$ is given by:
$$
2^{2n-1} n! \left(\xi_{2n}(ix)+\xi_{2n}(-ix)\right).
$$

\bigskip

One can further check that

$$
H_n(x)=i^n2^{n-1}n! \left((-1)^n\xi_n(-ix) +\xi_n(ix)\right).
$$

\bigskip

Now we come to the main point of this section.

\bigskip

It is well known that we are in the presence of a bispectral situation involving the differential operator 
 $\mathcal{H}_x$ and the difference operator  $\mathcal{B}_n$ given by

$$
(\mathcal{B}_n f)_n=f_{n+1}+2nf_{n-1},
$$

\bigskip

i.e., we have  
\begin{equation}\label{EqDifHermite}
\left(\mathcal{H}_xH_n\right)(x)=-2nH_n(x),
\end{equation}

\bigskip


as well as 
\begin{equation}\label{EqRecurrHermite}
\left(\mathcal{B}_nH_n\right)(n)=2xH_n(x).
\end{equation}

\bigskip

The subindices $x$ and $n$ going with the operators above are meant to emphasize that they act on the variables $x$ and $n$ respectively.

\bigskip

For this pair of operators 
$\mathcal{H}_x$
and  $\mathcal{B}_n$ more is actually true: for  
each $n=1,2,3,...$ one can find a {\bf nonpolynomial} eigenfunction $\psi_n(x)$ of $\mathcal{H}_x$ with eigenvalue $-2n$  such that

\begin{equation*}\label{recurrencia_Hermite}
\mathcal{B}_n \psi_n=\psi_{n+1}+2n\psi_{n-1}=2x\psi_n.
\end{equation*}
\bigskip

The family $\psi_n(x)$ with this recursion property is uniquely determined up to a global nonzero constant $b_0$ and one can see that
\begin{eqnarray*}
\psi_0(x)&=& \dfrac{\sqrt{\pi}\Erf (i x)}{2i}=x+\frac{x^3}{3}+\frac{x^5}{10}+\frac{x^7}{42}+\ldots,\\
\psi_1(x)&= &(-1)\left(1-x^2-\frac{x^4}{6}-\frac{x^6}{30}-\ldots \right), \\
\psi_2(x) &= &(-4)\left(x-\frac{x^3}{3}-\frac{x^5}{30}-\frac{x^7}{210}-\ldots \right). 
\end{eqnarray*}

\bigskip

Combining these results one can see that for any choice of constants $a_0,b_0$ the family of eigenfunctions
of $\mathcal{H}_x$
\begin{equation}\label{comblineal}
a_0H_n(x)+b_0\psi_n(x),\quad n\geq 1,
\end{equation}gives a bispectral situation going with the pair of operators 
$\mathcal{H}_x$ and $\mathcal{B}_n$.

\bigskip

The search for not necessarily polynomial solutions of a differential operator related to the Hermite one, starts in \cite{Reach123}. In the case of a differential operator going with one case of exceptional orthogonal Hermite polynomials, one can see in \cite{GLM} a certain recursion relation satisfied by a nonpolynomial sequence of solutions of the differential equation. This gives two different rank one cases and not a rank two case as we are looking for.

\bigskip

 The rank of the operators in a bispectral situation, such as the one defined by the continuous operator $\mathcal{H}_x$ in (\ref{EqDifHermite}) and the discrete operator $\mathcal{B}_n$ in (\ref{EqRecurrHermite}), is the dimension of the vector space of common eigenfunctions. In this particular case, for any fixed value of the discrete parameter $n$, $n \geq 1$, the resulting space of solutions in the continuous variable $x$ is two-dimensional, as given by (\ref{comblineal}). The straightforward definition given above is all that is needed here. However, we mention that this is connected with definitions of rank in the case of commutative rings of differential operators, starting with [6]. A very good account is in \cite{W1} (see also \cite{W2}) as well as in \cite{Haine1}, where Luc Haine reviews joint work of his with Plamen Iliev.
\bigskip

An intrinsic description of these $\psi_n(x)$ is possible: the functions $\psi_n(x)$ are odd or even in an alternating fashion, and since the Hermite polynomials are even or odd in the same complimentary
alternating fashion this characterizes the $\psi_n(x)$, since each function in $\mathcal{V}_n$ can be decomposed into the sum of an even and an odd function.

\bigskip
	
	Using the hypergeometric functions introduced above, the functions $\psi_n$ can be written as:
	\begin{eqnarray}
\psi_{2n+1}(x)&=&(-1)^{n+1}4^nn! \, _1F_1\left(-n-\frac{1}{2};\frac{1}{2};x^2\right), \label{laspsi_impar} \nonumber\\
\psi_{2n}(x)&=&	(-4)^n n! x\, _1F_1\left(\frac{1}{2}-n;\frac{3}{2};x^2\right), \quad n\geq0. \label{laspsi_par} \nonumber
	\end{eqnarray}


\bigskip

In terms of the function introduced earlier we have the expressions for $\psi_n(x)$

\[
\psi_{2n+1}(x)=(-1)^{n-1}2^{3n-1}n!\sqrt{\pi}(2n+1)!!\left[\xi_{2n+1}(ix)+\xi_{2n+1}(-ix)\right], 
\]
and
\[
\psi_{2n}(x)=i\sqrt{\pi}(-1)^n2^{3n-2}n!(2n-1)!!\left[\xi_{2n}(ix)-\xi_{2n}(-ix)\right], \quad n\geq 0.
\]


\bigskip


\bigskip

Once we have this bispectral family, then as we observed above we can get a collection of two dimensional families by picking $a_0$ and $b_0$ arbitrarily.

In the continuous-continuous case discussed in \cite{DG} the situation connected with the Korteweg-deVries equation features a pair of rank one bispectral families while the case related to the (so called) even family, later identified in \cite{ZM} to go with the master symmetries of KdV, goes with a full  rank two bispectral family.

\bigskip

In the language of \cite{DG,GH1}  we see that in the Hermite situation considered in this section we have
a family of rank two situations: one can pick an eigenfunction in $\mathcal{V}_0$ arbitrarely, by choosing the constants $a_0$ and $b_0$,
and this determines a family of eigenfunctions of (1) going with $\mathcal{B}_n$ and agreeing with it for $n=0$.

\bigskip

An important observation in \cite{DG} is that the rank one families are connected with what is called trivial monodromy while the case of a single rank two situation is connected to the case with the integers $\mathbb{Z}$ as the monodromy group. For a very recent look at the important issue of monodromy one can see \cite{FeldVa}, where the problem first addressed in \cite{DG} is taken a fresh look, and connected it with several other areas of mathematical physics.

\bigskip

The results in \cite{DG} cannot be taken as a general rule in terms of the relation between the Darboux process going from one differential operator to another and the corresponding monodromy groups. The peculiar situation in \cite{DG} is that a sequence of Darboux transformations preserves the mo\-nodromy group in two different instances, in one case the groups are all trivial and in the other case they coincide with the integers.
The paper \cite{G99} looks at an example of the hypergeometric equation of Euler, Gauss and Kummer with a finite monodromy group (all of these were characterized by H. Schwartz) and one sees that an application of the Darboux process to it gives a differential operator with an infinite monodromy group.

\bigskip

The result in \cite{DG} describing all cases with trivial monodromy was extended to the Hermite case in \cite{Oblom}. The example above (and those coming in later sections) is one of these and yet we have a single rank two bispectral family. The lesson is that one cannot assume that the situation in the case of the harmonic oscillator
 is the same as the case of the free particle considered in \cite{DG}.

\bigskip

We close this section with a few observations.
For a different and earlier look at the Hermite situation one can see \cite{McKTru}.
Very important papers related to the results in \cite{DG} are \cite{AdMos} and \cite{AiMcMos} as well as the older papers \cite{BurCha, Crum}. 


\bigskip

In the case of the Legendre operator one has the explicit expression for the
``Legendre functions of the second kind'' given below. It is clear that they satisfy the same recursion relation as the polynomials, since one has

\bigskip

$$
Q_n(y)=\dfrac{1}{2}\int_{-1}^1\dfrac{P_n(y)}{z-y}dy.
$$

\medskip

This is much simpler than the situation in the Hermite case, there is no need to do any scaling, as  in the expressions for
$\psi_{2n+1}$ and $\psi_{2n}$ above starting from the most familiar solutions of the differential equation, namely $P_n$ and $Q_n$.

\bigskip

Finally, our operator $\mathcal{H}_x$ can be conjugated by means of $$\dfrac{1}{h(x)}\mathcal{H}_x h(x)=\partial_x^2-x^2+1$$ in a more standard form, using $h(x)=e^{\tfrac{x^2}{2}}$. 
\section{Example 1: Hermite with one gap}

In this section we consider the operator
\begin{equation}\label{opdif_Honegap}
L=\partial_x^2-2(x+\frac{1}{x})\partial_x
\end{equation}
which is obtained from the original Hermite operator $\mathcal{H}_x$ by one application of the Darboux process, namely the factorization
\begin{equation*}\label{Darbouxonegap}
\mathcal{H}_x+2I=\partial_x^2-2x\partial_x+2=(\dfrac{1}{x}\partial_x-2)(x\partial_x-1)
\end{equation*}
gives rise to the operator $L$ defined by
\begin{eqnarray*}
L+2I&\equiv &(x\partial_x-1)(\dfrac{1}{x}\partial_x-2)\\
&=& \partial_x^2-2\left( x+\dfrac{1}{x}\right)\partial_x+2.
\end{eqnarray*}

It is easy to see that the polynomial eigenfunctions $p_n(x)$ of $L$, i.e., polynomial solutions of
$$
\left( L+2n \right)p_n(x)=0,
$$
can be written in the usual Wronskian form
$$
p_n(x)=\det\begin{pmatrix} x&1\\H_n(x)&H^{\prime}_n(x) \end{pmatrix}.
$$

\bigskip

We see that

\begin{eqnarray*}
p_0(x)&=&-1,\\
p_1(x)&=&0,\\
p_2(x)&=&4x^2+2,\\
p_3(x)&=&16x^3,\\
p_4(x)&=&48x^4-48x^2-12,\\
p_5(x)&=&128x^5-320x^3,\\
p_6(x)&=&120+720x^2-1440x^4+320x^6.
\end{eqnarray*}

It is easy to see, either by direct computation or by using the discrete-continuous version of the Wilson bispectral involution, originally given in \cite{W1}, and nicely put to use in \cite{BHY1,BHY2,KR,KM}, that one has the recursion
\begin{equation*}\label{recursion}
\dfrac{n(n-1)}{2}p_{n-2}(x)+\dfrac{2n-1}{4}p_n(x)+\dfrac{n-1}{8(n+1)}p_{n+2}(x)=\dfrac{x^2}{2}p_n(x),\quad n\geq 1,
\end{equation*}

\bigskip

assuming that $p_{-1}(x)=0$.

\bigskip

We observe that the fact that $\tau(x)=\left(\dfrac{x^2}{2} \right)^{\prime}=\Theta'(x)$, in the notation of \cite{DG} and \cite{Reach123}, vanishes at $x=0$ renders these polynomials as an instance of a bispectral family, but not quite a family of exceptional orthogonal polynomials. If these polynomials were ``exceptional orthogonal'' ones, then the orthogonality weight would be $e^{-x^2}$ divided by the square of this ``tau function''. The fact that this example is not one of orthogonal exceptional polynomials can be seen in \cite{GUGrMi}. More about this issue in the next section.

For a small sample of references on this topic see \cite{Dur-Herm,Dcorr,Du,GU1,GU2,GUGrMi,GUKKM,KM,Qu,STZ}. It may be worth remarking that the first use of the Darboux process in the context of exceptional orthogonal polynomials appears to be \cite{Qu} where the author very graciously says ``This procedure whose idea goes back to an older paper \cite{DG}", while the first observation of the presence of recursion relations appears to be \cite{STZ}.

\bigskip

In the spirit of the previous section we consider the full two dimensional space of solutions of
\begin{equation}\label{difeqphi}
\left( L\phi_n \right)(x)=-2n\phi_n(x)
\end{equation} and exhibit examples of non-polynomial families $\phi_n(x)$, indexed by $n$, that give rise to bispectral situations. The solutions to (\ref{difeqphi}) are given as follows
$$
\phi_n(x)=\alpha_n\left( 1+\sum_{i=1}^{\infty} S_{n,i} x^{2i} \right)+\beta_n \left( x^3+\sum_{i=2}^{\infty} \Gamma_{n,i} x^{2i+1} \right),
$$

where the scalars $\alpha_n$ and $\beta_n$ are arbitrary and the coefficients $S_i$ and $\Gamma_i$ are given as follows

\begin{equation*}
S_{n,i}=\left\{ \begin{array}{cc} n,&\text{if}\ i=1, \\\dfrac{(-1)^{i+1}}{(i)!(2i-3)!!}\displaystyle\prod_{k=0}^{i-1}(n-2k),&\text{if}\ i>1,\end{array}\right.
\end{equation*} and
\begin{equation*}
\Gamma_{n,i}=\dfrac{3(-1)^{i+1}}{(i-1)!(2i+1)!!}.\prod_{k=1}^{i-1}(n-2k-1).
\end{equation*} 
\bigskip

We are assuming that $\phi_{-2}=\phi_{-1}=0$.

\bigskip

One can write
the functions $\phi_n(x)$ in terms of the hypergeometric function, namely
$$
\phi_n(x)=\alpha_n\, _1F_1\left(-\frac{n}{2};-\frac{1}{2};x^2\right)+\beta_nx^3 \, _1F_1\left(\frac{3}{2}-\frac{n}{2};\frac{5}{2};x^2\right).
$$

\bigskip

In fact it is easy to see that
$$
(x\partial_x-1)_1F_1\left(-\frac{n}{2};\frac{1}{2};x^2\right)=-\ _1F_1\left(-\frac{n}{2};-\frac{1}{2};x^2\right)
$$

\bigskip

as well as

$$
(x\partial_x-1)\left[x\ _1F_1\left(-\frac{n}{2}+\frac{1}{2};\frac{3}{2};x^2\right)\right]=\dfrac{2}{3}(1-n)x^3\ _1F_1\left(\frac{3}{2}-\frac{n}{2};\frac{5}{2};x^2\right),
$$

\noindent highlighting the important role of the operator $x\partial_x-1$ that appears in the factorization
of $\mathcal{H}_x$ above, in giving new eigenfunctions in terms of old ones.


\bigskip

For the issue at hand, which is bispectrality beyond the family of polynomials given above, the question becomes: are there good choices of the sequences $\alpha_n$ and $\beta_n$ such that the functions $\phi_n(x)$ would satisfy the same recursion relation as the one satisfied by the polynomials $p_n(x)$?

\bigskip

We find that this is the case exactly when
$$
\alpha_n=\left\{ \begin{array}{cc} \widetilde{\alpha}_n,& n=0,1,\\ & \\-2(n-1)\alpha_{n-2}&n\geq 2, \end{array}\right.
$$ and $$
\beta_n=\left\{ \begin{array}{cc} \widetilde{\beta}_n,& n=0,\ 1,\ 2,\ 3,\\ &\\ 
-\dfrac{2(n-1)n}{n-3}\beta_{n-2}, &   n\geq 4. \end{array}\right.
$$



\bigskip


For illustration we display the first few $\phi_n(x)$:

\begin{eqnarray*}
	\phi_0(x)&=&\widetilde{\alpha}_0+\widetilde{\beta}_0\left( x^3+\dfrac{3}{5}x^5+\ldots \right),\\
	\phi_1(x)&=&\widetilde{\alpha}_1 \left(1+x^2+\dfrac{x^4}{2}+\dfrac{x^6}{6}+\dfrac{x^8}{24}+\ldots \right)
	+	\widetilde{\beta}_1\left( x^3+\dfrac{2}{5}x^5+\dfrac{4}{35}x^7+\dfrac{8}{315}x^9+\ldots \right),\\
	\phi_2(x)&=&-2\widetilde{\alpha}_0 \left(1+2x^2\right)+\widetilde{\beta}_2\left(x^3+\dfrac{x^5}{5}+\ldots\right),\\
		\phi_3(x)&=&-4\widetilde{\alpha}_1 \left(1+3x^2-\dfrac{3}{2}x^4+\ldots \right)+\widetilde{\beta}_3x^3,\\
			\phi_4(x)&=&12\widetilde{\alpha}_0 \left(1+4x^2-4x^4\right)-24\widetilde{\beta}_2\left(x^3-\dfrac{x^5}{5}+\ldots \right), \\
				\phi_5(x)&=&32\widetilde{\alpha}_1 \left(1+5 x^2+\ldots \right) -20\widetilde{\beta}_3\left(x^3-\dfrac{2}{5}x^5 \right), \\
				\phi_6(x)&=&-120\widetilde{\alpha}_0 \left(1+6x^2-12x^4+\dfrac{8}{3}x^6\right) +480\widetilde{\beta}_2\left(x^3-\dfrac{3}{5}x^5+\dfrac{3}{70}x^7+\ldots \right). \\
\end{eqnarray*}

	\bigskip

We notice that the functions $\phi_{i}(x)$ are the sum of a polynomial part and a nonpolynomial part, with the exception of $\phi_{1}$. This is connected to the fact that $p_{1}=0$. Naturally, all of this is related to the fact that the operator $L$ in (\ref{opdif_Honegap}) is obtained from $\mathcal{H}_x$ in (\ref{opHermite}) by the Darboux method. Notice that the polynomial parts coincide with $p_{j}(x)$ up to a scalar.

\bigskip

Notice that $\mathcal{H}_x$ as well as $L$ are invariant under the reflection that sends $x$ into $-x$. This played an important role in expressing the eigenfunctions of $L$  in terms of even or odd functions.

\bigskip
	
	Here, the coefficients $\widetilde{\alpha}_i$, $i=0,1$ and $\widetilde{\beta}_i$, $i=0,1,2,3$, are arbitrary as long as we want to satisfy the recursion relation starting with $n=4$.
If one wants the recursion to hold for smaller values of $n$ we need to impose restrictions on these parameters.

	
	
 If one wants to obtain the family of polynomials $p_n(x)$ given above  with these bispectral functions, one should choose
$$
\widetilde{\alpha}_0=-1,\ \widetilde{\alpha}_1=0, \ \widetilde{\beta}_0=0,\ \widetilde{\beta}_1=0,\ \widetilde{\beta}_2=0    \ \textrm{and} \ \widetilde{\beta}_3=16.
$$



\bigskip


One should point out that our expression for $L$ can be conjugated by means of the function $h(x)=x e^{x^2/2}$ into the operator


$$
\dfrac{1}{h(x)}L h(x)=\partial^2_{x}-x^2-1+2\dfrac{d^2}{dx^2}\log x.
$$

\bigskip

Regarding the question of the monodromy of this operator,
the classification of all those operators with quadratic grow in the potential was given, as mentioned earlier, in \cite{Oblom}. The author proved a conjecture of A. Veselov.
According to this classification the present case is one of trivial monodromy.

\bigskip

Once again, in spite of having trivial monodromy, we can exhibit a rank-two bispectral family, showing that the situation here differs from that in \cite{DG}.
In fact, there are several families of such rank-two bispectral objects parametrized by the parameters $\widetilde\alpha$ and $\widetilde\beta$ above.

\section{Example 2: Hermite with two gaps}

In this section we consider the operator
\begin{equation}\label{operator_H2gaps}
\widetilde{L}=\partial_x^2-\frac{2 \left(2 x^3+3 x\right)}{2 x^2-1}\partial_x+\dfrac{8}{2x^2-1},
\end{equation} which can be obtained from the original Hermite operator $\mathcal{H}_x$ in (\ref{opHermite}) by  application of the Darboux process, as will be seen below.

\bigskip

One can see that $$S \widetilde{L}= \mathcal{H}_{x} S,$$ if $S$ is given by

$$ S=  \dfrac{1}{2x^2-1}\partial_x^2  -\dfrac{4x}{2x^2-1}\left(1+\dfrac{1}{2x^2-1}\right)\partial_x\
+\dfrac{4}{2x^2-1}\left(1+\dfrac{1}{2x^2-1}\right)+2.$$

\bigskip

While this is true, it is not too helpful.

\bigskip

If we define $T$ as

$$T=(2x^2-1) \partial_x^2- 4x \partial_x,$$

we get the more useful relation

 $$ \widetilde{L} T = T \mathcal{H}_{x}.$$

\bigskip

This is the relation that says that $\widetilde{L}$ is obtained from $\mathcal{H}_x$ by an application of the Darboux method.

\bigskip

It is easy to see that the polynomial eigenfunctions $q_n(x)$ of $\widetilde{L}$, i.e., polynomial solutions of
$$
\left( \widetilde{L}+2n \right)q_n(x)=0,
$$
can be written by means of a determinant
$$
q_n(x)=\det\begin{pmatrix} 4x^2-2&8x\\H_{n-1}(x)&H^{\prime}_{n-1}(x) \end{pmatrix}.
$$

\bigskip

The first few $q_n(x)$ are given by

\begin{eqnarray*}
	q_0(x)&=&0,\\
	q_1(x)&=&-8x,\\
	q_2(x)&=&-8x^2-4,\\
	q_3(x)&=&0,\\
	q_4(x)&=&32x^4+24,\\
	q_5(x)&=&128x^5-128x^3+96x.
\end{eqnarray*}

\bigskip

The operator $T$ defined above converts the eigenfunctions of $\mathcal{H}_x$ into the
$q_n$ given above up to a scalar multiple. Of course this is one of the pay-offs of the Darboux process.

\bigskip

It is easy to see, either by direct computation or by using the tools based on the Wilson bispectral involution, see \cite{W1}, that one has the recursion
$$
\frac{2}{3} (n-1) (n-2) (n-3)q_{n-3}(x)+(n-3)(n-1)q_{n-1}(x)+\dfrac{(n - 3)}{2}q_{n+1}(x)+
\dfrac{(n - 3)}{12n}q_{n+3}(x)=\left(\dfrac{2x^3}{3}-x\right)q_n(x),
$$ which holds true for $n\geq 1$, assuming $q_{-1}(x)=0$ and $q_{-2}(x)=0$.

\bigskip

Observe that the corresponding ``tau function'' in this case is given by $\tau(x)= 2x^2-1$, which vanishes at two
real values of $x$. This gives an example of a family of Hermite polynomials with two gaps that cannot be
orthogonal exceptional ones. This should be compared with the example of Section $4$.  In other words, the condition in \cite{GUGrMi} about the even codimension in the Hermite case is neccesary but not sufficient.

\bigskip

Once again, in the spirit of the previous sections we consider the full two dimensional space of solutions of
\begin{equation}\label{difeqphi2}
\left( \widetilde{L}\phi_n \right)(x)=-2n\phi_n(x)
\end{equation}and exhibit examples of non-polynomial families $\phi_n(x)$, indexed by $n$, that give rise to bispectral situations. The solutions to (\ref{difeqphi2}) are given as follows
$$
\phi_n(x)=\alpha_n\left( 1+\sum_{i=1}^{\infty} S_{n,i} x^{2i} \right)+\beta_n \left( x+\sum_{i=1}^{\infty} \Gamma_{n,i} x^{2i+1} \right),\quad n\geq 0.
$$where the scalars $\alpha_n$ and $\beta_n$ are arbitrary and the coefficients $S_i$ and $\Gamma_i$ are given as follows

\begin{equation*}
S_{n,i}=\left\{ \begin{array}{cc} -(n-4),&\text{if}\ i=1, \\\dfrac{(-1)^{i}(n+4i(i-2))}{i!(2i-1)!!}\displaystyle\prod_{k=1}^{i-1}(n-2k),&\text{if}\ i>1,\end{array}\right.
\end{equation*} and
\begin{equation*}
\Gamma_{n,i}=\left\{ \begin{array}{cc} -(n-1)/3,&\text{if}\ i=1, \\ \dfrac{(-1)^{i}(n-1)(n+(2i-3)(2i+1))}{(2i+1)i!(2i-1)!!}\prod_{k=2}^{i-1}(n-2k-1),&\text{if}\ i>1.\end{array}\right.
\end{equation*} 
\bigskip



We are assuming that $\phi_{-2}=\phi_{-1}=0$.

\bigskip

	We can write the functions $\phi_n(x)$ in terms of the hypergeometric function, namely
	
	\begin{eqnarray*}
&&	\phi_n(x)=\alpha_n\, \left[\left(1-2 x^2\right) \, _1F_1\left(-\frac{n}{2};\frac{1}{2};x^2\right)+\left(6 x^2-4 x^4\right) \, _1F_1\left(1-\frac{n}{2};\frac{3}{2};x^2\right)\right]+\\
	&&\beta_n\left[ x \, _1F_1\left(\frac{1}{2}-\frac{n}{2};\frac{3}{2};x^2\right)+\frac{4}{15} (n-1) x^5 \, _1F_1\left(\frac{5}{2}-\frac{n}{2};\frac{7}{2};x^2\right) \right].
		\end{eqnarray*}
	

\bigskip

Once again, for the issue of bispectrality the question becomes: are there good choices of the sequences $\alpha_n$ and $\beta_n$ such that the functions $\phi_n(x)$ would satisfy the same recursion relation as the one satisfied by the polynomials $p_n(x)$?

\bigskip

We find that this requires, in an analogy to the previous example
$$
\alpha_n=\left\{ \begin{array}{cc} \widetilde{\alpha}_n,& n=0,1,2\\ & \\-2(n-1)\alpha_{n-2}&n\geq 3, \end{array}\right.
$$ and $$
\beta_n=\left\{ \begin{array}{cc} \widetilde{\beta}_n,& n=0,\ 1,\ 2,\ 3,\ 4,\ 5 \\ &\\ 
-\dfrac{2(n-3)(n-2)}{n-5}\beta_{n-2}, &   n\geq 6. \end{array}\right.
$$

\bigskip


For illustration we display the first few $\phi_n(x)$:

\begin{eqnarray*}
	\phi_0(x)&=&\widetilde{\alpha}_0(\left(  1+4x^2-\dfrac{16}{15}x^6-\dfrac{64}{105}x^8+\ldots\right)+\widetilde{\beta}_0\left( x+\dfrac{x^3}{3}-\dfrac{x^5}{6}-\dfrac{x^7}{6}-\dfrac{5}{72}x^{9}+\ldots \right),\\
	\phi_1(x)&=&\widetilde{\alpha}_1 \left(1+3x^2-\dfrac{x^4}{6}-\dfrac{13}{30}x^6-\dfrac{11}{56}x^8+\ldots \right)+\widetilde{\beta_1}x,\\
	\phi_2(x)&=&\widetilde{\alpha}_2 \left(1+2x^2\right)+\widetilde{\beta_2}\left(x-\dfrac{x^3}{3}+\dfrac{7}{30}x^5+\dfrac{23}{210}x^7+\dfrac{47}{1512}x^9+\ldots\right),\\
	\phi_3(x)&=&-4\widetilde{\alpha}_1 \left(1+x^2+\dfrac{x^4}{2}+\dfrac{x^6}{6}+\dfrac{x^8}{24}+\ldots \right)+\widetilde{\beta}_3\left(x-\dfrac{2}{3}x^3+\dfrac{8}{15}x^5+\dfrac{16}{105}x^7+\dfrac{32}{945}x^9+\ldots\right),\\
	\phi_4(x)&=&-6\widetilde{\alpha}_2 \left(1+\dfrac{4}{3}x^4\right)+\widetilde{\beta}_4\left(x-x^3+\dfrac{9}{10}x^5+\dfrac{5}{42}x^7+\dfrac{7}{360}x^9+\ldots \right), \\
	\phi_5(x)&=&32\widetilde{\alpha}_1 \left(1- x^2+\dfrac{5}{2}x^4-\dfrac{17}{30}x^6-\dfrac{37}{840}x^8+\ldots \right) +\widetilde{\beta}_5\left(x-\dfrac{4}{3}x^3+\dfrac{4}{3}x^5 \right). \\
\end{eqnarray*}

\bigskip

Notice that, once again, the functions $\phi_i$ have a different character than the remaining ones when $i=0$ and $i=3$, corresponding to indices when $p_i$ vanishes. For the polynomial parts we see that they are scalar multiples of the $q_i$, as in the previous example.

\bigskip

If we are satisfied with a recursion valid for $n$ starting at $6$ one has to impose 
	$$\widetilde{\beta}_{n+3}=-4n(n+1)\widetilde{\alpha}_n,\quad n=0,1,2,$$
	 and we have still lots of free parameters giving two dimensional bispectral families.  If we want the recursion to hold for smaller values of $n$ some further  restrictions have to be imposed among these parameters.

\bigskip

The operator $\widetilde{L}$ can be conjugated as follows
$$\dfrac{1}{h(x)}\widetilde{L}h(x)=\partial_x^2-x^2-3+2\dfrac{d^2}{dx^2} \log(2x^2-1),$$ using $h(x)=e^{\tfrac{x^2}{2}}(2x^2-1)$. This allows us to write our operator as a perturbation of the Hermite one using a Wronskian expression involving Hermite polynomials. In this case the polynomials are $H_0$ and $H_2$, in the previous example the polynomials are $H_0$ and $H_1$.



\section{The ad-conditions}

In this section we collect some interesting Lie-algebraic relations satisfied by each of the examples above. The arguments in \cite{DG}, properly adapted in \cite{Reach123}, show that they are intimately connected to bispectrality. In the continuous-continuous case of \cite{DG} they are necessary and sufficient, in the continuos-discrete case of \cite{Reach123} they are at least necessary and may be sufficient too.

We write
\begin{equation*}
\textrm{ad}\ X(Y)=XY-YX 
\end{equation*}for the usual commutator and \begin{equation*}\left(\textrm{ad}\ X\right)^n(Y)=\textrm{ad}\ X \left(\left(\textrm{ad}\ X\right)^{n-1}(Y)\right),\quad  n\geq 2.\end{equation*}

\bigskip

In Section 3, with 
$$
\Theta(x)=x
$$

and $\mathcal{L}$ the operator denoted by $\mathcal{H}_x$ in (\ref{opHermite}), we have the operator identity
$$
 \textrm{ad}\
\mathcal{L}^2(\Theta) -
4\textrm{ad}\
\mathcal{L}^0(\Theta)=0,
$$

which, except for simple changes can be found in \cite{Reach123}.

For the example in Section 4 with $\Theta(x)=\dfrac{x^2}{2}$ and $\mathcal{L}$ the operator $L$ in (\ref{opdif_Honegap}), we have the identity
$$ \textrm{ad}\ \mathcal{L}^3(\Theta) - 16\textrm{ad}\ \mathcal{L}^1(\Theta)=0.$$

\bigskip

One may suspect that this may be a consequence of a stronger identity
$$
\textrm{ad}\ \mathcal{L}^2(\Theta) - 16 \textrm{ad}\ \mathcal{L}^0(\Theta)=0,
$$

but this is not true. As usual $\textrm{ad}\ \mathcal{L}^0(\Theta)$ denotes the operator of multiplication by $\Theta$.

\bigskip

Finally, for the example of Section 5, choosing
$$
\Theta(x)=\dfrac{2}{3}x^3-x,
$$
and with $\mathcal{L}$ the operator $\widetilde{L}$ in (\ref{operator_H2gaps}), one gets the operator identity

$$
\textrm{ad}\ \mathcal{L}^4(\Theta) - 40 \textrm{ad}\ \mathcal{L}^2(\Theta)+144\textrm{ad}\ \mathcal{L}^0(\Theta)=0.
$$
\bigskip

We close with the observation that given one of the ad-conditions satisfied above, reversing the engine and trying to solve it for $\mathcal{L}$ and $\Theta$ is a formidable task. Solving the equations might give new examples. In the original case in \cite{DG}, solving explicitly  for $\mathcal{L}$ and $\Theta$ the equation $$\textrm{ad}^3\mathcal{L}\left(\Theta \right) =0$$was the key to led to the rational solutions of the Korteweg-deVries equation and its master symmetries.

\end{document}